\documentclass[a4paper, reqno]{amsart}%{{{1
\usepackage{amssymb}
\usepackage[utf8]{inputenc}
\usepackage[T1]{fontenc}
\usepackage{lmodern}
\usepackage{ae}
\usepackage[color=blue!30]{todonotes}
\usepackage{graphicx}
\usepackage{enumerate}
\usepackage{fancybox}
\usepackage{setspace}
\usepackage[all,xdvi]{xy}
\usepackage{ifthen}
\usepackage{everypage}
\usepackage{tikz}
\usepackage{hyperref}
\usetikzlibrary{shapes,matrix,arrows,calc}
\title{Virtual tangles and fiber functors}
\author{Adrien Brochier}\address{MPIM Bonn, Germany}\email{abrochier@mpim-bonn.mpg.de}
\date{\today}

\tikzstyle cross=[preaction={draw=white, -, line width=4pt}, thick]
\tikzstyle vcross=[preaction={draw=blue!30, -, line width=4pt}, thick]
\tikzstyle normal=[thick]
\tikzstyle chord=[densely dotted, thick]
\tikzstyle zero=[ultra thick, gray]
\tikzstyle zerocross=[preaction={draw=white, -, line width=4pt}, ultra thick, gray]
\tikzstyle point=[draw,circle,inner sep=1,fill=black]
\tikzstyle petitpoint=[draw,circle,inner sep=0.3,fill=black]

\newcommand{\negative}[3][-]{\draw[normal,#1] (#2+1,-#3).. controls (#2+1,-#3-0.3) and (#2,-#3-0.7)..(#2,-#3-1); \draw[cross,#1] (#2,-#3).. controls (#2,-#3-0.3) and (#2+1,-#3-0.7)..(#2+1,-#3-1);}
\newcommand{\positive}[3][-]{ \draw[normal,#1] (#2,-#3).. controls (#2,-#3-0.3) and (#2+1,-#3-0.7)..(#2+1,-#3-1);\draw[cross,#1] (#2+1,-#3).. controls (#2+1,-#3-0.3) and (#2,-#3-0.7)..(#2,-#3-1);}

\newcommand{\up}[3][-]{ 
\draw[normal,#1] (#2,-#3).. controls (#2,-#3-0.3) and (#2+0.2,-#3-0.5).. (#2+0.5,-#3-0.5);
\draw[normal] (#2+0.5,-#3-0.5).. controls (#2+0.8,-#3-0.5) and (#2+1,-#3-0.3).. (#2+1,-#3);}
\newcommand{\ap}[3][-]{
\draw[normal,#1] (#2,-#3-1).. controls (#2,-#3-0.7) and (#2+0.2,-#3-0.5).. (#2+0.5,-#3-0.5);
\draw[normal] (#2+0.5,-#3-0.5).. controls (#2+0.8,-#3-0.5) and (#2+1,-#3-0.7).. (#2+1,-#3-1);}

\newcommand{\straight}[3][-]{\draw[normal,#1] (#2,-#3) -- (#2,-#3-1);}

\newcommand{\vpositive}[3][-]{ 
\fill[blue!30] (#2-0.2,-#3)rectangle (#2+1.2,-#3-1);
\fill[white] (#2+0.5,-#3) circle (0.3);
\fill[white] (#2+0.5,-#3-1) circle (0.3);
\draw[normal,#1] (#2,-#3).. controls (#2,-#3-0.3) and (#2+1,-#3-0.7)..(#2+1,-#3-1);\draw[vcross,#1] (#2+1,-#3).. controls (#2+1,-#3-0.3) and (#2,-#3-0.7)..(#2,-#3-1);
}

\newcommand{\filledpositive}[3][-]{ 
\fill[blue!30] (#2-0.2,-#3)rectangle (#2+1.2,-#3-1);
\draw[normal,#1] (#2,-#3).. controls (#2,-#3-0.3) and (#2+1,-#3-0.7)..(#2+1,-#3-1);\draw[vcross,#1] (#2+1,-#3).. controls (#2+1,-#3-0.3) and (#2,-#3-0.7)..(#2,-#3-1);
}
\newcommand{\filling}[2]{
	\fill[blue!30] (#1+0.2,-#2) rectangle (#1+0.8,-#2-1) ;
}
\newcommand{\fillednegative}[3][-]{ 
\fill[blue!30] (#2-0.2,-#3)rectangle (#2+1.2,-#3-1);
\draw[normal,#1] (#2+1,-#3).. controls (#2+1,-#3-0.3) and (#2,-#3-0.7)..(#2,-#3-1); \draw[vcross,#1] (#2,-#3).. controls (#2,-#3-0.3) and (#2+1,-#3-0.7)..(#2+1,-#3-1);
}

\newcommand{\filledstraight}[3][-]{
\fill[blue!30] (#2-0.2,-#3)rectangle (#2+0.2,-#3-1);
\draw[normal,#1] (#2,-#3)--(#2,-#3-1) ;
}
\newcommand{\triup}[3][-]{ 
\fill[blue!30] (#2-0.2,-#3)rectangle (#2+1.2,-#3-1);
\begin{scope}
\clip (#2-0.2,-#3)rectangle (#2+1.2,-#3-1);
\fill[white] (#2+0.5,-#3) ellipse (0.3 and 0.6);
\end{scope}

\draw[normal] (#2,-#3) -- (#2,-#3-1);
\draw[normal] (#2+1,-#3) -- (#2+1,-#3-1);
}
\newcommand{\tridown}[3][-]{ 
\fill[blue!30] (#2-0.2,-#3)rectangle (#2+1.2,-#3-1);
\begin{scope}
\clip (#2-0.2,-#3)rectangle (#2+1.2,-#3-1);
\fill[white] (#2+0.5,-#3-1) ellipse (0.3 and 0.6);
\end{scope}
\draw[normal] (#2,-#3) -- (#2,-#3-1);
\draw[normal] (#2+1,-#3) -- (#2+1,-#3-1);
}
\newcommand{\virtual}[3][-]{ 
\draw[line width=.4cm, blue!30] (#2,-#3).. controls (#2,-#3-0.3) and (#2+1,-#3-0.7)..(#2+1,-#3-1);
\draw[normal,#1] (#2,-#3).. controls (#2,-#3-0.3) and (#2+1,-#3-0.7)..(#2+1,-#3-1);
\draw[preaction={draw=white, -, line width=0.5cm}, line width=.4cm, blue!30] (#2+1,-#3).. controls (#2+1,-#3-0.3) and (#2,-#3-0.7)..(#2,-#3-1);
\draw[normal,#1] (#2+1,-#3).. controls (#2+1,-#3-0.3) and (#2,-#3-0.7)..(#2,-#3-1);
}
\newcommand{\virtualop}[3][-]{
	\begin{scope}[xscale=-1]
		\virtual[#1]{#2}{#3}
	\end{scope}
}
\newcommand{\virtualup}[3][-]{ 
\draw[line width=.4cm, blue!30,#1] (#2,-#3).. controls (#2,-#3-0.3) and (#2+0.2,-#3-0.5).. (#2+0.5,-#3-0.5);
\draw[line width=.4cm, blue!30] (#2+0.5,-#3-0.5).. controls (#2+0.8,-#3-0.5) and (#2+1,-#3-0.3).. (#2+1,-#3);
\draw[normal,#1] (#2,-#3).. controls (#2,-#3-0.3) and (#2+0.2,-#3-0.5).. (#2+0.5,-#3-0.5);
\draw[normal] (#2+0.5,-#3-0.5).. controls (#2+0.8,-#3-0.5) and (#2+1,-#3-0.3).. (#2+1,-#3);}
\newcommand{\virtualap}[3][-]{
\draw[line width=.4cm, blue!30,,#1] (#2,-#3-1).. controls (#2,-#3-0.7) and (#2+0.2,-#3-0.5).. (#2+0.5,-#3-0.5);
\draw[line width=.4cm, blue!30] (#2+0.5,-#3-0.5).. controls (#2+0.8,-#3-0.5) and (#2+1,-#3-0.7).. (#2+1,-#3-1);
\draw[normal,#1] (#2,-#3-1).. controls (#2,-#3-0.7) and (#2+0.2,-#3-0.5).. (#2+0.5,-#3-0.5);
\draw[normal] (#2+0.5,-#3-0.5).. controls (#2+0.8,-#3-0.5) and (#2+1,-#3-0.7).. (#2+1,-#3-1);}

\newcommand{\oldvirtual}[3][-]{ 
	\draw[normal,#1] (#2,-#3).. controls (#2,-#3-0.3) and (#2+1,-#3-0.7)..(#2+1,-#3-1);\draw[normal,#1] (#2+1,-#3).. controls (#2+1,-#3-0.3) and (#2,-#3-0.7)..(#2,-#3-1);\draw (#2+0.5,-#3-0.5) circle (0.1);
}
\newcommand{\filledup}[3][-]{ 
\fill[blue!30] (#2-0.2,-#3)rectangle (#2+1.2,-#3-1);
\draw[normal,#1] (#2,-#3).. controls (#2,-#3-0.3) and (#2+0.2,-#3-0.5).. (#2+0.5,-#3-0.5);
\draw[normal] (#2+0.5,-#3-0.5).. controls (#2+0.8,-#3-0.5) and (#2+1,-#3-0.3).. (#2+1,-#3);}
\newcommand{\filledap}[3][-]{
\fill[blue!30] (#2-0.2,-#3)rectangle (#2+1.2,-#3-1);
\draw[normal,#1] (#2,-#3-1).. controls (#2,-#3-0.7) and (#2+0.2,-#3-0.5).. (#2+0.5,-#3-0.5);
\draw[normal] (#2+0.5,-#3-0.5).. controls (#2+0.8,-#3-0.5) and (#2+1,-#3-0.7).. (#2+1,-#3-1);}

\newcommand{\counit}[2]{
	\begin{scope}
	\clip (#1-0.2,-#2) rectangle (#1+1.2,-#2-1);
	\fill[blue!30] (#1+0.5,-#2) circle (0.7);
	\end{scope}
}
\newcommand{\unit}[2]{
	\begin{scope}
		\clip (#1-0.2,-#2) rectangle (#1+1.2,-#2+1);
	\fill[blue!30] (#1+0.5,-#2) circle (0.7);
	\end{scope}
}

\newcommand{\tik}[1]{\begin{tikzpicture}[baseline=(current bounding box.center)] #1 \end{tikzpicture} }
\numberwithin{equation}{section}

\newcommand{\RR}{\mathbb{R}}

\newcommand{\mf}[1]{\mathfrak{#1}}

\newcommand{\id}{\operatorname{id}}

\newcommand{\h}{\hslash}

\newcommand{\ot}{\otimes}

\newcommand{\T}{\mathcal T}
\newcommand{\B}{\mathcal B}

\newcommand{\cC}{\mathcal {C}}

\newcommand{\cV}{\mathcal {V}}
\newcommand{\cS}{\mathcal {S}}

\newcommand{\kk}{\mathbb K}

\newcommand{\eno}{\operatorname{End}}
\newcommand{\vect}{\operatorname{vect}}

\newcommand{\modu}{\operatorname{-mod}}

\newtheorem{thm}{Theorem}[section]
\newtheorem*{thm*}{Theorem}

\newtheorem{prop}[thm]{Proposition}
\newtheorem*{prop*}{Proposition}

\let\oldcite\cite
\renewcommand{\cite}[2][noOption]{\ifthenelse{\equal{#1}{noOption}}{{\scriptsize \oldcite{#2}}}{{\scriptsize \oldcite[#1]{#2}}}}
\newtheorem{defi}[thm]{Definition}

\theoremstyle{remark} 
\newtheorem{rmk}[thm]{Remark}

\begin{document}
\maketitle
\begin{abstract}
	We define a category $v\T$ of tangles diagrams drawn on surfaces with boundaries. On the one hand we show that there is a natural functor from the category of virtual tangles to $v\T$ which induces an equivalence of categories. On the other hand, we show that $v\T$ is universal among ribbon categories equipped with a strong monoidal functor to a symmetric monoidal category. This is a generalization of the Shum--Reshetikhin--Turaev theorem characterizing the category of ordinary tangles as the free ribbon category. This gives a straightforward proof that all quantum invariants of links extends to framed oriented virtual links. This also provides a clear explanation of the relation between virtual tangles and Etingof--Kazhdan formalism suggested by Bar--Natan. We prove a similar statement for virtual braids, and discuss the relation between our category and knotted trivalent graphs.
\end{abstract}
\setcounter{section}{-1}
\section{Introduction}%{{{
Braids and framed oriented tangles can be interpreted as morphisms in certain categories $\B$ and $\T$ respectively. These categories can be characterized by the following universal property:
\begin{thm*}[Shum~\cite{Shum1994}, Reshetikhin--Turaev~\cite{Reshetikhin1990}]
	The categories $\B$ and $\T$ are respectively the free braided monoidal category and the free ribbon category on one object.
\end{thm*}
We stress the fact that we do not assume the underlying monoidal category to be strict, so $\B$ and $\T$ are really the parenthesised versions of the categories of braids and tangles as introduced in~\cite{Bar-Natan1998,Bar-Natan1997}. 

The main result of this paper is a similar description of an appropriate category of \emph{virtual} tangles. 
The notion of \emph{virtual} knotted objects is due to Kauffman\cite{Kauffm1999}. It is a diagrammatic generalization of usual knotted objects in the following sense: usual knots can be represented by knot diagrams, which are \emph{planar} tetravalent graphs whose vertices are positive and negative crossings, modulo the Reidemeister relations. Virtual knot diagrams are simply obtained by dropping the planarity condition. In other words, the edges of the graph are allowed to intersect. These intersections are precisely the virtual crossings, usually depicted as follow:
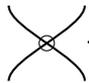
\[
	\tik{\oldvirtual{0}{0}}.
\]

One can instead choose an abstract punctured surface on which the non planar graph representing a given virtual knot can be drawn without self-intersection~\cite{CisnerosdelaCruz2015,Kuperberg2003,ScottCarter2002}. This point of view leads naturally to a topological interpretation of virtual knotted object, as link or tangle diagrams drawn on surfaces modulo homeomorphisms and certain additional relations (tearing and puncturing) which reflects the fact that the choice of the surface is not unique. This leads to the following alternative diagrammatic description of ordinary and virtual crossings
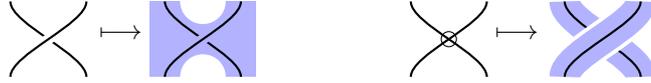
\begin{align*}
	 \tik{\positive{0}{0}} & \longmapsto\tik{\vpositive{0}{0}}&\tik{\oldvirtual{0}{0}} &\longmapsto \tik{\virtual{0}{0}}  
\end{align*}
where the blue strip is a piece of the underlying surface.

The main goal of this paper is to find a categorical interpretation of these pictures, extending the one existing for usual tangles. As for the original result, it provides both a natural framework for constructing virtual links invariants and a clear connection with deformation-quantization, namely with the quantum Yang--Baxter equation and Etingof--Kazhdan formalism.

It turns out that the answer is quite simple: in any graphical calculus for a certain kind of categorical structure, objects are represented by points while morphisms  are depicted by string diagrams. On the other hand, \emph{functors} between categories are depicted by intervals, and the functoriality on morphisms is expressed by embedding string diagrams into 2-dimensional strips (see e.g.~\cite{McCurdy2012}). Roughly speaking, the main observation of this paper is the fact that the pictures above are part of a graphical calculus for monoidal functors from a ribbon category to a symmetric monoidal category.

Hence, we introduce a certain category $v\T$ whose objects are finite sequences of non-associative words, and whose morphisms are punctured framed surfaces with tangles diagrams drawn on them modulo certain relations. We first show:
\begin{prop*}
	The category $v\T$ is a strict symmetric monoidal category. There is a natural functor from the category of parenthesized tangles $\T$ to $v\T$ given by ``embedding a tangle into a blue rectangle''. This functor carries a canonical strong monoidal structure. 
\end{prop*}
We then prove a Reidemeister-like theorem for those objects, leading to a presentation of $v\T$ by generators and relations. This implies our main result (see Theorem~\ref{thm:main}):
\begin{thm*}
	The category $v\T$ is universal for the properties of the above Proposition.
\end{thm*}
It means that the data of a ribbon category $\cC$, an object $V\in \cC$, a strict symmetric monoidal category $\cS$ and a strong monoidal functor from $\cC$ to $\cS$ gives rise to representations of $v\T$ extending the representations of $\T$ attached to $(\cC,V)$, i.e. to a functor $G_V:v\T\rightarrow \cS$ fitting in the following commutative diagram:
\begin{center}\begin{tikzpicture}[description/.style={fill=white,inner sep=2pt}] 
 \matrix (m) [matrix of math nodes, row sep=3em, 
 column sep=2.5em, text height=1.5ex, text depth=0.25ex] 
 {\T & v\T \\ 
\cC & \cS\\}; 
 \path[->,font=\scriptsize] 
 (m-1-1) edge node[auto] {$\iota$}(m-1-2) 
 (m-1-1) edge node[auto] {$F_V$}(m-2-1) 
 (m-1-2) edge node[auto] {$G_V$}(m-2-2) 
 (m-2-1) edge node[auto] {$G$}(m-2-2) 
 ; 
\end{tikzpicture}\end{center}

Let $v\T_0$ be the category whose objects are associative words on $\{+,-\}$ and morphisms are equivalence classes of framed oriented virtual tangles (see Definition~\ref{def:vtzero} for a precise definition), then the above topological interpretation induces a functor $v\T_0 \subset v\T$. We show that this is an equivalence of categories (Proposition~\ref{prop:vt0}). We stress the fact that this is, however, \emph{not} an isomorphism of categories. This is because $v\T$ contains morphisms of the form
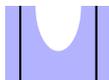
\[
	\tik{\triup{0}{0}}
\]
which are not virtual tangles, yet plays an important role in our story since they are precisely the structure maps of the monoidal structure on the functor $\T\rightarrow v\T$. This relates to the basic observation that the ordinary crossings once interpreted as diagrams on a surface as above are not elementary and can be decomposed into smaller pieces. We note that there is a \emph{strict} monoidal functor from the strict version $\T^{str}$ of the category of tangles to $v\T_0$. It is well known that $\T$ and $\T^{str}$ are equivalent as ribbon categories, and formally the above universal property could be equivalently stated for the pair $(\T^{str},v\T_0)$ instead. Yet, this would make its proof much less natural and would completely hide its relation with the graphical calculus of functors. Also, some of the most important examples of ribbon categories arising ``in nature'' are non-strict and it turns out that many results are much more natural and transparent when stated in terms of $\T$. In some sense one goal of this paper is to advocate $v\T$ as a better replacement of $v\T_0$ which account for a certain non-strictness which naturally arises in examples. 

In fact, additional morphisms in $v\T$ have a natural interpretation as virtual knotted trivalent graph (vKTG). It also gives a clean interpretation of the relation between usual long knotted trivalent graphs and parenthesised tangles and of some natural operations defined on vKTG's. This is discussed in Section~\ref{sec:KTG}.

The Shum--Reshetikhin--Turaev theorem is the basic ingredient tying quantum algebra and low-dimensional topology. In one direction, it shows that ribbon categories provide well-behaved links invariants, i.e. invariants which can be computed by cutting links into elementary pieces, and which are compatible with certain natural operations. The main examples of explicit ribbon categories comes from representation theory of quantum groups, leading to vast generalizations of the Jones polynomial. Our result implies in particular that links invariants derived from quantum groups extends naturally to virtual links (Proposition~\ref{prop:quantum-invariants}). It also provides a graphical calculus for those categories, allowing to do algebraic computation by drawing pictures.

In the other direction, this result is crucial in explaining the topological background of deformation-quantization: roughly speaking, certain analytically-defined topological invariants, inspired by conformal and quantum field theories, can be used to produce quantization of interesting structures. The main example of this situation is the Drinfeld-Kontsevich invariant~\cite{Drinfeld1990a,Kontsevich1993} (see also~\cite{Bar-Natan1997,Cartier1993,Kassel1998}): this is a functor from the category of tangles to a certain category of Feynman diagrams (so-called chord diagrams) constructed by formally integrating a universal version of the KZ connection. This functor therefore induces a ribbon structure on the category of diagrams; such structures are in one-to-one correspondence with Drinfeld associators. This construction is the key technical point of the proof of the formality of the little 2-discs operad which plays a prominent role in the modern approach to deformation-quantization. 

Representations of the category of chord diagrams diagrams can be obtained from the category of modules over a finite-dimensional Lie algebra $\mf g$ equipped with an element $t \in S^2(\mf g)^{\mf g}$. The Drinfeld--Kontsevich invariant then induces a highly non-trivial ribbon structure on the category $\mf g\modu[[\h]]$. There is a precise sense in which $t$ turns the symmetric monoidal category of $\mf g$-modules into some sort of Poisson category (see e.g.~\cite{Pantev2011}) of which this procedure gives a canonical quantization. 

A fundamental theorem of Drinfeld~\cite{Drinfeld1990a} asserts that if $\mf g$ is simple, then the ribbon category obtained this way is equivalent as a ribbon category to the category of modules over the corresponding quantum group. It gives a conceptual explanation of why the invariants attached to quantum groups exist but does not quite recover the quantum groups themselves. This question has been solved by Etingof--Kazhdan: any pair of a Lie algebra $\mf g$ and a solution $r$ of the classical Yang-Baxter equation leads to a pair $(\mf g,t)$ as above, hence to a ribbon structure on $\mf g\modu[[\h]]$. The main result of~\cite{Etingof1996} is that this ribbon category can be realized as the category of modules over a quasi-triangular Hopf algebra quantizing the pair $(\mf g,r)$. Then, if $\mf g$ is simple there is a standard choice for $r$ whose quantization recovers the quantum group attached to $\mf g$. The main technical point of the proof is the construction of a certain strong monoidal structure on the forgetful functor from the ribbon category $\mf g[[\h]]\modu$ to the category of vector spaces. This is ``Bar--Natan's dream''~\cite{Bar-Natan2011} that this construction should come from an analog of the Drinfeld-Kontsevich invariant for virtual tangles. Our result is one step in that direction.

%}}}
%}}}
\section{The category $v\T$}%{{{
\subsection{The tangle category}%{{{
Recall the following definitions from~\cite{Bar-Natan1997} to which we refer for details.
\begin{defi}
	A parenthesised word on a set $S$ is an element of the free non-associative monoid generated by $S$.
\end{defi}
If $w$ is a parenthesised word, denote by $\bar w$ the word obtained by forgetting the parenthesis.
\begin{defi}
	Let $\T$ be the category whose objects are parenthesised words on the set $\{+,-\}$, and morphisms from $w$ to $w'$ are isotopy classes of framed, oriented tangles whose source is $\bar w$ and target is $\bar w'$. Let $\B$ be the category whose objects are parenthesised words on $\{\bullet\}$ and morphisms are isotopy classes of braids.
\end{defi}
\begin{thm}\label{thm:SRT}
	The category $\T$ is a ribbon category: the tensor product is given by 
	\[
		(w,w')\mapsto ww',
	\]
	the associativity constraint by the identity tangle from $(w_1w_2)w_3$ to $w_1(w_2w_3)$, the braiding by the obvious braid diagram from $ww'$ to $w'w$. Objects $+$ and $-$ are dual to each other with structure maps given by the cup and the cap diagram. This category is the universal ribbon category in the following sense: for any ribbon category $\cC$ and any object $V\in \cC$, there exists a unique strict monoidal functor
	\[
F_V:\T\longrightarrow \cC
	\]
preserving the braiding and the ribbon structure, and such that $F_V(+)=V$ and $F_V(-)=V^*$.

The category $\B$ is a braided monoidal category, the braiding being defined similarly, and for every pair of a braided monoidal category $\cC$ and an object $V\in \cC$ there is a unique strict braided monoidal functor $F:\B\rightarrow \cC$ such that $F(\bullet)=V$.
\end{thm}
%}}}
\subsection{The virtual knotted objects category}%{{{
\begin{defi}\label{defi:vt}
	Let $v\T$ be the category defined as follow: objects are finite sequences of parenthesised words on $\{+,-\}$, which we denote by words enclosed in square bracket to avoid confusion with the concatenation:
\[
	[w_1]\dots [w_k].
\]

Morphisms between two objects
\begin{align*}
	W&=[w_1]\dots[w_k] & W'=[w_1']\dots [w_l']
\end{align*}
are non necessarily connected, framed compact surfaces with at least one circle boundary component on each connected component, together with:
\begin{itemize}
	\item a choice of a marked open interval on the boundary for each $w_i$ and each $w_i'$. We impose that near each interval marked by some $w_i$ (resp. some $w_i'$) the surface is locally framed-diffeomorphic to $\RR\times \RR_{\leq 0}$ (resp. $\RR\times \RR_{\geq 0}$) equipped with its standard framing.
	\item a framed oriented tangle diagram drawn on the surface, whose endpoints are attached to the marked intervals in a way consistent with the word written on the said interval in the obvious sense.
\end{itemize}
Here ``drawn on the surface'' means embedded in a small thickening of the surface. Morphisms are considered up to:
\begin{itemize}
	\item homeomorphisms of the underlying surface which preserves the boundaries, the framing and the diagrams
	\item isotopies of framed oriented tangles
	\item addition or removal of a connected component without any marked interval or diagram
	\item ``tearing'' and ``puncturing'' as explained below.
\end{itemize}
Composition is given by gluing matching intervals and tangles endpoints.
\end{defi}

Tearing means cutting the surface along a line joining two boundary components which does not intersect the tangle diagram. Conversely, given two intervals distinct from the marked intervals, on two non necessarily distinct boundary components of the surface we allow those intervals to be glued together, if it can be done in a way compatible with the framing. 
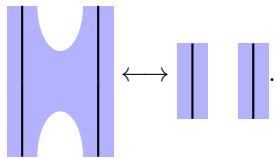
\[
\tik{\triup{0}{0}\tridown{0}{1}}\longleftrightarrow\tik{\filledstraight{0}{0}\filledstraight{0.8}{0}}.
\]
Puncturing means removing a disc which does not contains nor intersect a diagram. Conversely, we allow any non distinguished boundary component to be sealed by gluing a disc, if it can be done in a way compatible with the framing.
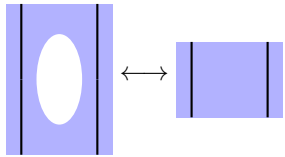
\[
\tik{\tridown{0}{0}\triup{0}{1}}\longleftrightarrow\tik{\filledstraight{0}{0}\filling{0}{0}\filledstraight{1}{0}}
\]

\begin{figure}
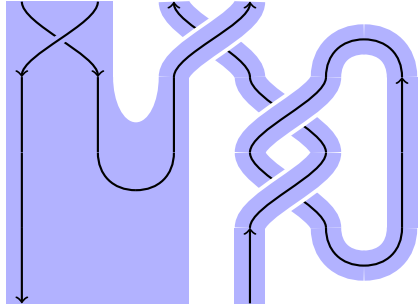
\label{fig:morph}
	\tik{
		\filledpositive[->]{0}{0} \virtual[<-]{2}{0}\virtualap{4}{0}
		\filledstraight{0}{1}\filling{0}{1}\triup{1}{1}\virtual{3}{1}\filledstraight[<-]{5}{1}
		\filledstraight{0}{2}\filling{0}{2}\filledup{1}{2}\virtual{3}{2}\filledstraight{5}{2}
		\filledstraight[->]{0}{3}\filling{0}{3}\filling{0.6}{3}\filling{1}{3}\filling{1.4}{3}\filledstraight[<-]{3}{3}\virtualup{4}{3}
	}
	\caption{A morphism from $[++][-][-]$ to $[+][-]$}
\end{figure}

A typical morphism in $v\T$ is shown on Figure~\ref{fig:morph}.

\begin{rmk}\label{rmk:convention}
When we draw such a morphism, it is assumed that the framing of the surface is the blackboard framing, and that the framing of a tangle diagram is the one induced by the framing of the surface, that the marked interval associated with $W$ and $W'$  appears on the top and the bottom of the pictures respectively, in the correct order. The direction of a morphism is from top to bottom. To clarifies pictures, we will denote by a single strand whose endpoints are labelled by the same word $w$, the identity diagram joining $w$ with itself in the unique way consistent with the orientation. For example, if $w=++-$ then:
\[
	\tik{
		\filledstraight{0}{0}
		\node[label=above:$w$] at (0,0) {};
		\node[label=below:$w$] at (0,-1) {};
	}:=
		\tik{\fill[blue!30] (-0.2,0)rectangle (2+0.2,-1);
	\draw[normal,->] (0,0)--(0,-1);
	\draw[normal,->] (1,0)--(1,-1);
	\draw[normal,<-] (2,0)--(2,-1);
}
\]
\end{rmk}

\section{A universal property for $v\T$}%{{{
In this section we define a symmetric strict monoidal structure on $v\T$ and construct a strong monoidal functor
\[
\T\longrightarrow v\T.
\]
We then show the main result of this paper, namely that $v\T$ is universal for these properties.

The tensor product on $v\T$ is given both on objects and morphisms by juxtaposition, and the unit is the empty word $\emptyset$\footnote{Not to be confused with $[\emptyset]$: these are isomorphic but different objects.}. The symmetric braiding is by successive applications of the virtual crossing, i.e. it is defined on the tensor product $[w_1]\ot[w_2]$ of two intervals by
\[
	\sigma_{[w_1],[w_2]}=\tik{\virtual{0}{0}
		\node[label=above:$w_1$] at (0,0) {};
		\node[label=above:$w_2$] at (1,0) {};
		\node[label=below:$w_2$] at (0,-1) {};
		\node[label=below:$w_1$] at (1,-1) {};
	}
\]
(using the conventions of Remark~\ref{rmk:convention}) and then defined iteratively in such a way that the hexagon axioms hold, e.g.
\[
	\sigma_{[w_1],[w_2][w_3]}=\tik{
		\virtual{0}{0}\filledstraight{2}{0}
		\filledstraight{0}{1}\virtual{1}{1}
		\node[label=above:$w_1$] at (0,0) {};
		\node[label=above:$w_2$] at (1,0) {};
		\node[label=above:$w_3$] at (2,0) {};
		\node[label=below:$w_2$] at (0,-2) {};
		\node[label=below:$w_3$] at (1,-2) {};
		\node[label=below:$w_1$] at (2,-2) {};
	}.
\]
\begin{prop}\label{prop:struct}
	Equipped with the above defined tensor product and symmetry, $v\T$ is a strict symmetric monoidal category.
\end{prop}
\begin{proof}
	Clear.
\end{proof}

As explained in the introduction, we want to relate the surfaces entering the definition of $v\T$ with the calculus of functor. This can be made precise as follow:
\begin{prop}\label{thm:functor}
There is a functor $\iota$ from $\T$ to $v\T$ defined on objects by
\[
	w\longmapsto [w]
\]
and on morphisms by embedding a tangle diagram into a blue rectangle. This functor is strong monoidal, with monoidal structure given by
\begin{equation}\label{eq:monoidal}
	\tik{\tridown{0}{0}}
\end{equation}
\end{prop}
\begin{proof}
The fact that $\iota$ is indeed a functor, i.e. that it is well defined and compatible with composition, is true by the very definition of $\T$ and $v\T$. Now again by definition the above picture is a morphism in $v\T$
\[
	\iota(w\ot w')\rightarrow \iota(w)\ot \iota(w').
\]
The fact that it is natural is clear by sliding pieces of tangles diagrams on the picture. This also imply that it is compatible with tensor products of morphisms since the latter is just given by concatenation.

This morphism is invertible, with inverse given by
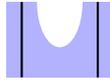
\[
	\tik{\triup{0}{0}}.
\]
Indeed, the puncturing relation implies that
\[
	\tik{\tridown{0}{0}\triup{0}{1}}=\tik{\filledstraight{0}{0}\filling{0}{0}\filledstraight{1}{0}}
\]
and the tearing relation implies

\[
	\tik{\triup{0}{0}\tridown{0}{1}}=\tik{\filledstraight{0}{0}\filledstraight{0.8}{0}}.
\]
The monoidal structure correspond to the following equality:
\[
	\tik{
	\filledstraight{0}{0}\filling{0}{0}\filledstraight{1}{0}\filling{1}{0}\filledstraight{2}{0}
		\tridown{0}{1}\filling{1}{1}\filledstraight{2}{1} 
		\filledstraight{0}{2}\tridown{1}{2}

		\node at (0.8,-1) {$($};
		\node at (2.2,-1) {$)$};
		\node at (-0.2,-0) {$($};
		\node at (1.2,-0) {$)$};
	}
	=
\tik{
	\filledstraight{0}{0}\filling{0}{0}\tridown{1}{0}
	\tridown{0}{1}\filledstraight{1}{1}\filledstraight{2}{1}
		
		\node at (-0.2,0) {$($};
		\node at (1.2,0) {$)$};
	}.
\]
Finally, the following pictures
\[
	\tik{\unit{0}{0}},\tik{\counit{0}{0}}
\]
are the unit $\emptyset\rightarrow [\emptyset]$ and the counit $[\emptyset]\rightarrow \emptyset$. The defining relation of $v\T$ implies that
\[
	\tik{\counit{0}{0}\unit{0}{1.5}}=\tik{\filling{0}{0}}=\id_{[\emptyset]}
\]
and
\[
	\tik{\unit{0}{0}\counit{0}{0}}=\tik{\draw[dotted] (0,0) rectangle (1,1);} =\id_{\emptyset}.
\]

\end{proof}
%}}}

The main result of this paper is essentially that Proposition~\ref{prop:struct} and~\ref{thm:functor} characterize the category $v\T$. 

\begin{thm}\label{thm:main}
	Let $\cC$ be a ribbon category and $V \in \cC$. Let $\cS$ be a strict symmetric monoidal category and 
	\[
G:\cC\longrightarrow \cS
	\]
	a strong monoidal functor. Then there exists a unique strict symmetric monoidal functor
	\[
G_V:v\T\longrightarrow \cS
	\]
	such that $G_V([+])=G(V)$ and such that the following diagram of functor is strictly commutative:
\begin{center}\begin{tikzpicture}[description/.style={fill=white,inner sep=2pt}] 
 \matrix (m) [matrix of math nodes, row sep=3em, 
 column sep=2.5em, text height=1.5ex, text depth=0.25ex] 
 {\T & v\T \\ 
\cC & \cS\\}; 
 \path[->,font=\scriptsize] 
 (m-1-1) edge node[auto] {$\iota$}(m-1-2) 
 (m-1-1) edge node[auto] {$F_V$}(m-2-1) 
 (m-1-2) edge node[auto] {$G_V$}(m-2-2) 
 (m-2-1) edge node[auto] {$G$}(m-2-2) 
 ; 
\end{tikzpicture}\end{center}
\end{thm}
\begin{proof}
Assuming that $G_V$ is well defined, the above conditions determines it completely, e.g.
\[
	G_V([+-])=G(V\ot V^*)
\]
and
\[
	G_V([+][-])=G(V)\ot G(V^*).
\]
Hence we need to show that $G_V$ is well defined. It follows from the following presentation of $v\T$:
\begin{thm}
Morphisms in $v\T$ are generated under composition and tensor product by the following elementary morphism:
\begin{align}
	\tik{\filledpositive{0}{0}},\tik{\fillednegative{0}{0}} , \tik{\filledup{0}{0}},\tik{\filledap{0}{0}}\label{gen:braid}\\
	\tik{\virtual{0}{0}},\tik{\virtualop{0}{0}}\label{gen:sym}\\
	\tik{\counit{0}{0}},\tik{\unit{0}{0}},	\tik{\triup{0}{0}},\tik{\tridown{0}{0}}\label{gen:mono}
\end{align}
and the following relations:
\begin{itemize}
	\item relations stating that the generators~\eqref{gen:sym} induce a symmetric monoidal structure, i.e. the fact that they are natural in both variable and that they satisfies the hexagon axioms together with
		\[
	\tik{\virtual{0}{0}}=\tik{\virtualop{0}{0}}\
		\]
	\item relations stating that $\iota$ is a well defined strong monoidal functor, i.e. the axioms of ribbon categories for generators~\eqref{gen:braid} except that we require naturality only with respect to morphisms which are images of morphisms in $\T$ through $\iota$, and relations appearing in the proof of Proposition~\ref{thm:functor}.
\end{itemize}
\end{thm}
\begin{rmk}
	Here we again follow the convention explained in Remark~\ref{rmk:convention}, i.e. in addition to all the generators above, we consider all morphisms obtained from those by either erasing or doubling a strand and with all possible orientations. 
\end{rmk}
\begin{proof}
	First of all observe that we did not include the following morphisms as generators:
	\[
		\tik{\virtualup{0}{0}}\ \tik{\virtualap{0}{0}}
	\]
	because they can be obtained by composing other generators, e.g.
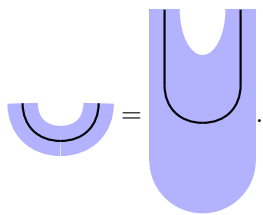
	\[
		\tik{\virtualup{0}{0}}=\tik{\triup{0}{0}\filledup{0}{1}\counit{0}{2}}.
	\]
	Including the morphisms above, that this set is generating is clear: forgetting about the tangles diagrams, every oriented surface with boundaries and marked intervals is homeomorphic to the thickening of a uni-trivalent fat graph (the univalent vertices accounting for the marked intervals), and adding curls if necessary one can arrange so that the framing match the blackboard framing and so that the top and bottom marked interval are at the right position (i.e. as in Figure~\ref{fig:morph}). This gives a way of writing every morphism in $v\T$ as a composition of the elementary morphisms above. Note that putting the marked interval at the correct position may force a part of the surface to overlap another but this can be done using~\eqref{gen:sym}. The tearing and puncturing relation and the removal of connected components without diagrams follows from the monoidal structure on $\iota$. The equality
	\[
		\tik{\triup{0}{0} \filling{1}{0}\filledstraight{2}{0}
		\filledstraight{0}{1} \filling{0}{1} \tridown{1}{1}}	
		=
		\tik{\filledstraight{0}{0}\tridown{1}{0}
		\triup{0}{1} \filledstraight{2}{1}	}	
	\]
	follows from the well-known fact that strong monoidal functors are in particular Frobenius monoidal functors. The fact that diagrams can slide on the underlying surface is either the functoriality of $\iota$, the fact that it is also monoidal with respect to tensor products of morphisms or the naturality of the symmetric braiding depending on the situation. Then the theorem follows from a combination of Theorem~\ref{thm:SRT} and the well known classification of surfaces.

\end{proof}
\end{proof}
\begin{rmk}\label{rmk:vb}
	One can define a variant $v\B$ as the category whose objects are finite sequences of parenthesised words on $\{\bullet\}$ and morphisms are defined as in Definition~\ref{defi:vt} but with oriented surfaces and braid diagrams instead (what we mean by ``braid diagram'' here can be made precise along the same line as in~\cite{CisnerosdelaCruz2015}). Then the analog of Theorem~\ref{thm:main} holds with $\B$ instead of $\T$ and for $\cC$ a braided monoidal category.
	\end{rmk}
%}}}
\section{Virtual tangles and virtual knotted trivalent graphs}%{{{
\subsection{Virtual tangles}%{{{
Tangles can be understood combinatorially as a certain set of \emph{planar} uni-tetravalent graphs modulo the Reidemeister relations. Virtual tangles are defined the same way but without the planarity condition: the edges of the graph can intersect, and those intersections are called virtual crossings. In the previous sections we used the fact that tangles, more precisely parenthesised tangles diagrams, form a category. Here we define a similar structure for virtual tangles. Since morphisms in categories are represented by planar graphs we need, even if it is slightly less natural, to consider the virtual crossing as an extra generator and add extra relations.
\begin{defi}\label{def:vtzero}
	Let $v\T_0$ be the \emph{strict} monoidal category whose objects are associative words on $\{+,-\}$ and morphisms are compositions of the following diagrams (with all possible orientation and the usual rule of compatibility between the orientation and the objects attached to the endpoints):
	\[
		\tik{\positive{0}{0}},\tik{\negative{0}{0}},\tik{\oldvirtual{0}{0}},\tik{\ap{0}{0}},\tik{\up{0}{0}}
	\]
	modulo the following relations
	\begin{itemize}
		\item planar isotopy
		\item the Reidemeister moves for framed oriented tangles
		\item the virtual Reidemeister moves
			\begin{align*}
				\tik{\oldvirtual{0}{0}\oldvirtual{0}{1}} &= \tik{\straight{0}{0} \straight{1}{0}} &
				\tik{\oldvirtual{0}{0} \straight{2}{0}
					\straight{0}{1} \oldvirtual{1}{1}
				\oldvirtual{0}{2} \straight{2}{2}}
				&=
				\tik{\straight{0}{0} \oldvirtual{1}{0}
					\oldvirtual{0}{1} \straight{2}{1}
				\straight{0}{2} \oldvirtual{1}{2}}&
							\end{align*}
		\item the mixed relation
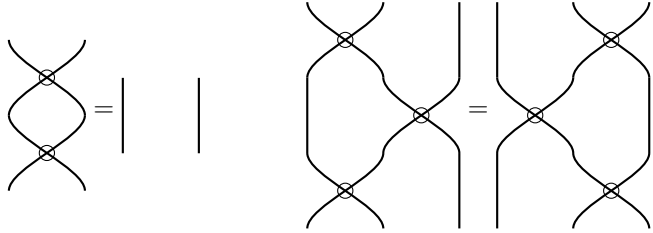
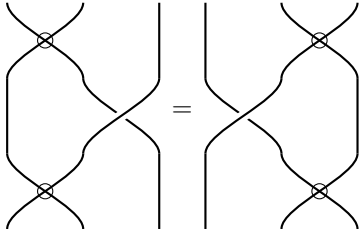
			\[
	\tik{\oldvirtual{0}{0} \straight{2}{0}
					\straight{0}{1} \positive{1}{1}
				\oldvirtual{0}{2} \straight{2}{2}}
				=
				\tik{\straight{0}{0} \oldvirtual{1}{0}
					\positive{0}{1} \straight{2}{1}
				\straight{0}{2} \oldvirtual{1}{2}}
			\]
	\end{itemize}
\end{defi}
%}}}
\begin{rmk}
	Note that we do not mod out by the virtual analog of the first Reidemeister move. This notion is called ``rotational virtual tangle'' in~\cite{Kauffm1999}.
\end{rmk}
\begin{prop}\label{prop:vt0}
The category $v\T_0$ has a natural symmetric monoidal structure. The functor $v\T_0 \rightarrow v\T$ defined on objects by
\[
	\epsilon_1\epsilon_2\dots \epsilon_n\longmapsto [\epsilon_1][\epsilon_2]\dots [\epsilon_n]
\]
and on morphisms by
\begin{align}
	 \tik{\positive{0}{0}} & \longmapsto\tik{\vpositive{0}{0}}&\tik{\oldvirtual{0}{0}} &\longmapsto \tik{\virtual{0}{0}} \\
	\tik{\ap{0}{0}} & \longmapsto \tik{\virtualap{0}{0}} &\tik{\up{0}{0}} & \longmapsto \tik{\virtualup{0}{0} }
\end{align}
is strict monoidal and induces an equivalence of categories $v\T_0 \simeq v\T$.
\end{prop}
\begin{proof}
	The symmetric monoidal structure is once again induced by concatenation and the virtual crossing.
	The fact that this functor is well defined i.e. that the defining relations of $v\T_0$ are satisfied is easily checked with pictures. Let $\cV$ be the full sub-category of $v\T$ whose objects are sequences of intervals containing exactly one letter. We need to show that this functor is surjective on morphisms. Given a morphism in $\cV$ one can detach every handle which does not contain a diagram thanks to the tearing relation, then remove every empty component. View the tangle diagram at hand as a graph whose vertices are either positive and negative crossings or endpoints of an open component. Then by tearing and puncturing if necessary one can arrange so that the underlying surface is a deformation-retract of this graph. Note that here we use in a crucial way the fact each marked interval contains exactly one endpoint of a tangle diagram: the previous statement is clearly not true for an arbitrary morphism in $v\T$. This shows that the morphism we started with is the image of a virtual tangle, i.e. that the morphisms in $\cV$ are generated under composition and tensor product by the pictures images of the generators of virtual tangles as shown in the Proposition. 
	
	Using Theorem~\ref{thm:main} one can show that the defining relations of $\cV$ are equivalent to the virtual Reidemeister moves, hence the functor is fully faithful. On the other hand, by repeated use of
	\[
		\tik{\triup{0}{0}},\tik{\tridown{0}{0}}
	\]
every object in $v\T$ is isomorphic to an object in $\cV$, so that this functor is essentially surjective.
\end{proof}
\begin{rmk}
	There is an obvious analog of this Proposition leading to an equivalence from the category of virtual braids and the category $v\B$ of Remark~\ref{rmk:vb}.
\end{rmk}
\subsection{Knotted trivalent graphs}\label{sec:KTG}%{{{
The goal of this section is to relate morphisms in $v\T$ with certain virtual knotted trivalent graphs (KTG)~\cite{Bar-Natan2014,Bar-Natan2013,Thurston2002}. This comparison will be mostly informal, since KTG's are described in the language of planar algebras, and unlike what happen for tangles the author is not aware of an appropriate categorical definition of those, or rather this section should be thought of as a proposal for such a definition.

The construction is similar to the proof of Proposition~\ref{prop:vt0} but starting with an arbitrary morphism in $v\T$. By tearing and puncturing appropriately, every such morphism can be seen as a virtual tangle, composed on the top and the bottom by a (possibly empty) succession of tensor products of
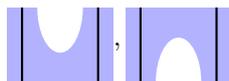
\[
\tik{\triup{0}{0}},\tik{\tridown{0}{0}}
\]
and identities, in an order dictated by the parenthesization of the words on the marked intervals. Then one reverse the map of Proposition~\ref{prop:vt0} on the middle tangle, and replace the successive compositions of the above pieces by trees.

 Here is an illustration of this process ; the first equality is by puncturing. Now the top-left piece of the figure, say, is the composition of the monoidal structure and its inverse. The first one turns into a trivalent vertex, while the second one is absorbed when ``reversing'' the map of Proposition~\ref{prop:vt0}. The trees are determined by the parenthesization of the top and bottom object.
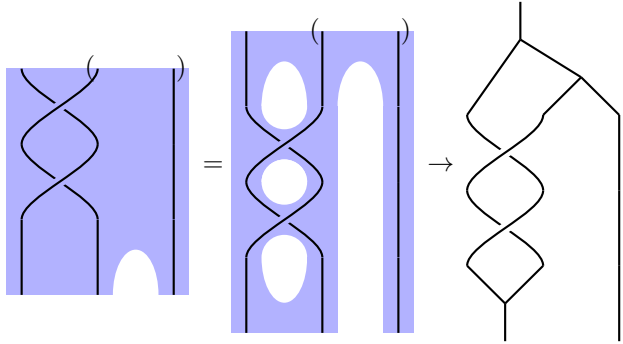
\[
	\tik{
		\filledpositive{0}{0}\filling{1}{0}\filledstraight{2}{0}
		\filledpositive{0}{1}\filling{1}{1}\filledstraight{2}{1}
		\filledstraight{0}{2}\filling{0}{2}\tridown{1}{2}
		\node at (0.9,0) {$($};
		\node at (2.1,0) {$)$};
}= 
	\tik{
		\tridown{0}{-1}\tridown{1}{-1}
		\vpositive{0}{0}\filledstraight{2}{0}
		\vpositive{0}{1}\filledstraight{2}{1}
		\triup{0}{2}\filledstraight{1}{2}\filledstraight{2}{2}
		\node at (0.9,1) {$($};
		\node at (2.1,1) {$)$};
		}
	\rightarrow
	\tik{
		\positive{0}{0}\straight{2}{0}
		\positive{0}{1}\straight{2}{1}
		\straight{2}{2}
		\draw[normal] (1.5,0.5)--(1,0);
		\draw[normal] (1.5,0.5)--(2,0);
		\draw[normal] (0.7,1)--(0,0);
		\draw[normal] (0.7,1)--(1.5,0.5);
		\draw[normal] (0.5,-2.5)--(0,-2);
		\draw[normal] (0.5,-2.5)--(1,-2);
		\draw[normal] (0.5,-2.5) -- (0.5,-3);
		\draw[normal] (0.7,1.5) -- (0.7,1);

	}
\]
Note that we could allow inner trivalent vertices as well, but they can all either be simplified or pushed on the top or the bottom of the diagram. Under this identification, composition corresponds to gluing the endpoint and applying the unzip operation of~\cite{Thurston2002} which in our setting corresponds to the tearing. We stress the fact that in our setting tearing gives equal morphisms, while in~\cite{Thurston2002} this is a transformation between isomorphic but different spaces of KTG's on different skeletons which cannot be directly compared. In the same way, puncturing corresponds to ``bubble cancellation'' as in~\cite{Bar-Natan2014}. Hence, we claim that every morphism space in $v\T$ can be identified with a space of equivalence classes of virtual KTG's under the unzip and bubble cancellation moves.

Applying this construction to the image of parenthesised tangles under $\iota$, this recovers Bar-Natan's observation that they are in bijection with long KTG's, i.e. KTG's without virtual crossings and only one end point at the top and at the bottom.
%}}}
%}}}
\section{Applications}%{{{

\subsection{Quantum invariants for virtual links}%{{{
In this section we explain how the formalism developed in the previous section allows one to extend quantum invariants of links to virtual links. Namely, let $H$ be a ribbon Hopf algebra. The main example, of course, is the quantum group $U_q(\mf g)$ associated to a semi-simple Lie algebra $\mf g$. The following is well known:
\begin{prop}
The category $H\modu$ of finite dimensional $H$-modules is a ribbon category. The forgetful functor which to an $H$-modules attaches the underlying vector space is strong monoidal, with monoidal structure the identity.
\end{prop}
Therefore, an immediate application of Theorem~\ref{thm:main} implies:
\begin{prop}\label{prop:quantum-invariants}
Each finite dimensional $H$-module $V$ gives rise to a numerical invariant for framed oriented virtual link whose restriction to classical links coincide with the Reshetikhin--Turaev invariant attached to $V$.
\end{prop}
\begin{rmk}
	Note that Theorem~\ref{thm:main} also implies that this constructions gives an $\eno(V)$-valued invariant of long virtual links, but this will \emph{not} be an $H$-module endomorphism in general, hence it cannot be identified with a scalar even if $V$ is simple.
\end{rmk}
\begin{rmk}
As for the classical case, the categorical formalism implies that those invariants are compatible with composition. In particular the above virtual links invariant is multiplicative on disjoint unions.
\end{rmk}

%}}}
\subsection{Etingof--Kazhdan formalism}%{{{
Let $(\mf g,r)$ be a finite-dimensional quasi-triangular Lie bialgebra. This means that $r\in \mf g^{\ot 2}$ is such that $t=r+r^{2,1}$ is $\mf g$-invariant and such that $r$ satisfies the classical Yang-Baxter equation
\[
	[r^{1,2},r^{1,3}]+[r^{1,2},r^{2,3}]+[r^{1,3},r^{2,3}]=0.
\]

Let $U(\mf g)[[\h]]\modu$ be the symmetric monoidal category of topologically free $U(\mf g)[[\h]]$-modules of finite type. Together with $t$, the choice of a Drinfeld associator  $\Phi$ turns this category into a ribbon category $U(\mf g)[[\h]]\modu_{\Phi}$. If one choose the so-called KZ associator~\cite{Drinfeld1990a} this recovers the Drinfeld--Kontsevich universal finite type invariant.

Let $F:U(\mf g)[[\h]]\modu_{\Phi}\rightarrow \vect_\h$ be the forgetful functor which sends a module to its underlying $\kk[\h]]$-module. The is a natural algebra isomorphism $\eno(F)\cong U(\mf g)[[\h]]$. One of the main result of~\cite{Etingof1996} is the construction, for each choice of $\Phi$, of a strong monoidal structure $J$ on $F$ which satisfies $J=1+\h r+O(\h^2)$. 

Therefore Theorem~\ref{thm:main} implies the following topological interpretation of the element $J$:
\begin{prop}
	Let $(\mf g, r)$ be a finite dimensional quasi-triangular Lie algebra. Each choice of a Drinfeld associator $\Phi$ and a $\mf g$-module $V$ gives rise to a functor
	\[
v\T\rightarrow \vect_\h
	\]
	which extends the Drinfeld-Kontsevich invariant attached to the same data.
\end{prop}

%}}}
%}}}
%}}}

\end{document}